\documentclass{amsart}

\usepackage{amssymb,amsthm,amsmath,amsxtra}
\usepackage{multicol}
\usepackage[all]{xy}

\theoremstyle{plain}
\newtheorem{thm}{Theorem}[section]
\newtheorem{prop}[thm]{Proposition}

\newtheorem{lem}[thm]{Lemma}
\newtheorem{cor}[thm]{Corollary}
\newtheorem*{thm*}{Theorem}
\newtheorem*{prop*}{Proposition}

\theoremstyle{remark}

\newtheorem{rmk}[thm]{Remark}

\newcommand{\psmod}[1]{~(\textup{\text{mod}}~{#1})}

\newcommand{\T}{\rule{0pt}{2.6ex}}
\newcommand{\Aprim}{A_{\text{prim}}}

\newenvironment{enumalph}
{\begin{enumerate}}
{\end{enumerate}}

\newcommand{\F}{\mathbb F}

\newcommand{\Q}{\mathbb Q}
\newcommand{\R}{\mathbb R}
\newcommand{\Z}{\mathbb Z}

\newcommand{\fraka}{\mathfrak a}
\newcommand{\frakd}{\mathfrak d}
\newcommand{\frakD}{\mathfrak D}
\newcommand{\frakf}{\mathfrak f}

\newcommand{\frakH}{\mathfrak H}
\newcommand{\frakM}{\mathfrak M}
\newcommand{\frakN}{\mathfrak N}
\newcommand{\frakp}{\mathfrak p}

\newcommand{\calD}{\mathcal D}
\newcommand{\calO}{\mathcal O}

\newcommand{\la}{\langle}
\newcommand{\ra}{\rangle}

\newcommand{\legen}[2]{\left(\frac{#1}{#2}\right)}

\DeclareMathOperator{\area}{area}

\DeclareMathOperator{\GL}{GL}
\DeclareMathOperator{\disc}{disc}

\DeclareMathOperator{\img}{img}

\DeclareMathOperator{\N}{N}
\DeclareMathOperator{\nrd}{nrd}
\DeclareMathOperator{\ord}{ord}
\DeclareMathOperator{\PSL}{PSL}

\DeclareMathOperator{\SL}{SL}

\title{Shimura curves of genus at most two}
\author{John Voight}

\begin{document}

\begin{abstract}
We enumerate all Shimura curves $X^\frakD_0(\frakN)$ of genus at most two: there are exactly $858$ such curves, up to equivalence.
\end{abstract}

\maketitle

The elliptic modular curve $X_0(N)$ is the quotient of the completed upper half-plane $\frakH^*$ by the congruence subgroup $\Gamma_0(N)$ of matrices in $\SL_2(\Z)$ that are upper triangular modulo $N \in \Z_{>0}$.  The curve $X_0(N)$ forms a coarse moduli space for (generalized) elliptic curves equipped with a cyclic subgroup of order $N$.  Modular curves have been studied in great detail on account of their importance in many fields, especially arithmetic geometry and the theory of automorphic forms.  The genus of $X_0(N)$ goes to infinity with $N$ according to an explicit formula \cite[Chapter~3]{DH}, \cite{CWZ}: for example, $X_0(N)$ has genus at most two if and only if $N \leq 29$ or $N \in \{31,32,36,37,49,50\}$.  Equations for modular curves of small genus have been explicitly determined in many cases \cite{MShimura}.  

In this paper, we consider these results in the context of Shimura curves, which are generalizations of modular curves.  A \emph{Shimura curve} $X_0^\frakD(\frakN)$ is the quotient of the (completed) upper half-plane $\frakH$ by the congruence arithmetic Fuchsian group $\Gamma_0^\frakD(\frakN) \subset \PSL_2(\R)$, which is specified by two coprime ideals of the ring of integers $\Z_F$ of a totally real number field $F$: the \emph{discriminant} $\frakD \subset \Z_F$, a squarefree ideal such that $n=[F:\Q]$ and the number of prime factors of $\frakD$ have opposite parity, and the \emph{level} $\frakN \subset \Z_F$ (see \S 1).  The curves $X_0^\frakD(\frakN)$, like modular curves, can be seen as moduli spaces for certain classes of abelian varieties equipped with level structure.  Taking $F=\Q$, $\frakD=\Z$, and $\frakN=N\Z$, one recovers the modular curves as Shimura curves.  

In this article, we consider Shimura curves of small genus.  In this vein, Long, Machlachlan, and Reid have proven the following theorem.

\begin{thm*}[Long-Maclachlan-Reid \cite{LMR}]
There are only finitely many congruence arithmetic Fuchsian groups $\Gamma$ of bounded genus, up to conjugacy.
\end{thm*}  

Already proving this statement in the simplest case---that there are only finitely many congruence subgroups of $\PSL_2(\Z)$ of genus $0$---is a nontrivial result of Thompson \cite{Thompson} (see also Chua-Lang-Yang \cite{CLY}).  The proof of Long et al.\ combines lower and upper bounds due to Selberg and Zograf, respectively, for the first nonzero eigenvalue of the Laplacian of $\Gamma$ (see \S 1 for precise statements), and thereby bounds the coarea of $\Gamma$ in terms of the genus.

The restriction to congruence subgroups is necessary: indeed, by the theory of Belyi covers, every curve defined over a number field is already uniformized by some subgroup of $\PSL_2(\Z)$, since $\PSL_2(\Z)$ contains $\Gamma(2)$, a free group on two generators.  By comparison, however, Takeuchi \cite{Takeuchi} has proven that there are only finitely many arithmetic Fuchsian groups with any fixed signature, and for certain simple signatures all such groups have been identified (see e.g.\ Machlachlan-Rosenberger \cite{MR}).  

Long et al.\ \cite{LMR} enumerate all maximal congruence arithmetic Fuchsian groups of genus zero with $F=\Q$ and produce a partial list of other such groups with $2 \leq [F:\Q] \leq 7$---however, they do not provide a complete list for any genus $g$.  Other partial lists have also been computed: Alsina-Bayer \cite{AB} list the invariants for Shimura curves corresponding to ``small ramified cases'' in the case $F=\Q$; and Johansson \cite{Johansson} has enumerated all Shimura curves of genus $g \leq 2$ with either (1) $F=\Q$ or (2) $[F:\Q]=2$ and $\frakN=\Z_F$.  As evidenced by these tables, it is reasonable to expect that the complete list of congruence arithmetic Fuchsian groups will be quite long already in the case of genus $g=0$. 

In this article, we consider the problem from the perspective of arithmetic geometry and enumerate all groups $\Gamma_0^\frakD(\frakN)$ corresponding to Shimura curves $X_0^\frakD(\frakN)$ of genus at most $2$.  We consider two such curves, specified by the ideals $\frakD,\frakN$ of $\Z_F$ and $\frakD',\frakN'$ of $\Z_{F'}$, to be \emph{equivalent} if there is an isomorphism $\sigma:F \xrightarrow{\sim} F'$ of fields such that $\sigma(\frakD)=\frakD'$ and $\sigma(\frakN)=\frakN'$.  

Our main result is as follows.

\begin{thm*}
There are exactly $858$ Shimura curves $X_0^\frakD(\frakN)$ of genus at most $2$, up to equivalence.
\end{thm*}

There are $258$ curves of genus $0$, $334$ of genus $1$, and $266$ of genus 2.  The complete list of these curves is provided in Tables 4.1--4.7 in Section 4.  If one fixes a unique field in each isomorphism class but applies no equivalence, then there are $1361$ such curves: $335,588,438$ of genus $0,1,2$, respectively.

The layout of the paper is as follows.  In \S 1, we apply the result of Selberg-Zograf to bound the hyperbolic area of a Shimura curve $X=X_0^\frakD(\frakN)$ of genus at most two; we then use the Odlyzko bounds to bound the root discriminant $\delta_F$ and thereby the degree $n=[F:\Q]$ of the totally real field $F$ in the data specifying $X$.  In \S 2, we discuss computational aspects of embedding numbers of orders in quaternion orders, quantities that appear in the genus formula for $X$.  Then in \S 3, for each possible field $F$, we enumerate the possibilities for the discriminant $\frakD$ and level $\frakN$, thereby enumerating all Shimura curves $X_0^\frakD(\frakN)$ of genus zero.  Finally, in \S 4 we tabulate the curves and in \S 5 we comment on some points of interest in the data.

The author wishes to thank Peter Clark and Victor Rotger for their encouragement in writing up this result, the Magma group at the University of Sydney for their hospitality, St\'ephane Louboutin for useful discussions, and the reviewer for his or her helpful comments.

\section{Notation and background}

In this section, we introduce the basic notions and background used throughout; basic references are given by Vign\'eras \cite{Vigneras}, Katok \cite{Katok}, and Elkies \cite{ElkiesSCC}.

Let $\overline{\Q}$ be an algebraic closure of $\Q$.  Let $F \subset \overline{\Q}$ be a totally real field of degree $[F:\Q]=n$ and let $d_F$ and $\delta_F=d_F^{1/n}$ denote the discriminant and root discriminant of $F$, respectively.  Let $B$ be a \emph{quaternion algebra} over $F$, a central simple $F$-algebra of dimension $4$.  There is a unique (anti-)involution $\overline{\phantom{x}}:B \to B$ such that the \emph{reduced norm} $\nrd(\gamma)=\gamma\overline{\gamma} \in F$ for all $\gamma \in B$.  The algebra $B$ is \emph{split} or \emph{ramified} at a place $v$ of $F$ according as $B \otimes_F F_v \cong M_2(F_v)$ or not, where $F_v$ denotes the completion of $F$ at $v$.  Let $S$ the set of places $v$ of $F$ where $B$ is ramified.  Then $S$ is a finite set of even cardinality, and determines $B$ uniquely up to $F$-algebra isomorphism.  Let $\frakD=\prod_{\frakp \nmid \infty,\ \frakp \in S} \frakp$ denote the \emph{discriminant} of $B$.  Suppose that $B$ is split at a unique real place, corresponding to $\iota_\infty:B \hookrightarrow B \otimes_F \R \cong M_2(\R)$.  

Let $\Z_F$ be the ring of integers of $F$.  An \emph{order} of $B$ is a full $\Z_F$-submodule of $B$ that is also a subring; an order is \emph{maximal} if it is not properly contained in any other.  Let $\calO \subset B$ be a maximal order.  Let $\frakN \subset \Z_F$ be an ideal coprime to the discriminant $\frakD$.  Choose an embedding 
\[ \iota_\frakN:\calO \hookrightarrow \calO \otimes_{\Z_F} \Z_{F,\frakN} \cong M_2(\Z_{F,\frakN}), \] 
where $\Z_{F,\frakN}$ denotes the completion of $\Z_F$ at $\frakN$.  The order
\[ \calO_0(\frakN) = \{\gamma \in \calO: \iota_\frakN(\gamma)\text{ is upper triangular modulo $\frakN$}\} \]
is known as an \emph{Eichler order} of level $\frakN$; it can also be obtained as the intersection of two maximal orders in $B$.  Let 
\[ \calO_0(\frakN)_1^*=\{\gamma \in \calO_0(\frakN) : \nrd(\gamma)=1\} \]
denote the group of units of $\calO(\frakN)$ of reduced norm $1$ and define
\[ \Gamma=\Gamma_0^\frakD(\frakN)=\{\iota_\infty(\gamma) : \gamma \in \calO_0(\frakN)_1^*\}/\{\pm 1\} \subset \PSL_2(\R) \]
and similarly $\calO(\frakN)=\{ \gamma \in \calO : \gamma \equiv 1 \pmod \frakN\}$ and 
\[ \Gamma^{\frakD}(\frakN)=\left\{\iota_\infty(\gamma):\gamma \in \calO(\frakN)_1^*\right\}/\{\pm 1\}. \]
We abbreviate $\Gamma_0^\frakD(\Z_F)=\Gamma_0^\frakD(1)$.  The quotient $X=X_0^\frakD(\frakN)=\Gamma_0^\frakD(\frakN) \backslash \frakH$ is a \emph{Shimura curve} and can be given the structure of a Riemann surface in a way that, up to isomorphism, does not depend on the choice of $\calO$, $\iota_\infty$, or $\iota_\frakN$.  

We have the formula of Shimizu \cite[Appendix]{Shimizu}:
\begin{equation} \label{shimizu}
A = \area(X)=\frac{4}{(2\pi)^{2n}} d_F^{3/2} \zeta_F(2) \Phi(\frakD) \Psi(\frakN),
\end{equation}
where $\zeta_F(s)$ denotes the Dedekind zeta function of $F$, and
\begin{eqnarray*}
 \Phi(\frakD)=\#(\Z_F/\frakD)^*=\N(\frakD)\prod_{\frakp \mid \frakD} \left(1-\frac{1}{\N(\frakp)} \right), \\
 \Psi(\frakN)=[\Gamma^\frakD(1):\Gamma_0^\frakD(\frakN)]=\N(\frakN)\prod_{\frakp \mid \frakN} \left(1+\frac{1}{\N(\frakp)} \right);
\end{eqnarray*}
here the hyperbolic area $\mu(D)=\displaystyle{\frac{1}{2\pi}\int\!\!\int_D \frac{dx\,dy}{y^2}}$ is normalized so that an ideal triangle has area $1/2$.  The genus $g$ of $X$ is determined by the Riemann-Hurwitz formula
\begin{equation} \label{RH}
A = 2g-2+\sum_q e_q\left(1-\frac{1}{q}\right)+e_\infty
\end{equation}
where $e_q$ is the number of elliptic cycles of order $q \in \Z_{\geq 2}$ and $e_\infty$ the number of parabolic cycles in $\Gamma$.  We accordingly say $\Gamma$ has \emph{signature} $(g;m_1,\dots,m_t;s)$ if $\Gamma$ has $t$ elliptic cycles of orders $m_1,\dots,m_t$ and $s$ parabolic cycles, and within a signature we abbreviate $m^k = \underbrace{m,\dots,m}_k$.  Note that $s=0$ unless $B=M_2(\Q)$, corresponding to the classical case of elliptic modular curves; if $s=0$, we write simply $(g;m_1,\dots,m_t)$.

We now state the Selberg-Zograf bound, following Long-Maclachlan-Reid \cite{LMR}.  A Fuchsian group $\Gamma$ is \emph{congruence} if $\Gamma \supset \Gamma^{\frakD}(\frakN)$ for some $\frakN$ and choice of $\calO$.  Let $\Gamma$ be a congruence arithmetic Fuchsian group of genus $g$ and area $A$, and let $\lambda_1(\Gamma)$ denote the first nonzero eigenvalue of the Laplacian of $\Gamma$.  Selberg \cite{Selberg} proved that the lower bound 
\[ \lambda_1(\Gamma) \geq 3/16 \] 
holds when $\Gamma$ is a congruence subgroup of $\SL_2(\Z)$; Vign\'eras \cite{Vigneras1-4}, relying on work of Jacquet-Langlands \cite{JL}, then showed how to generalize this result to all congruence arithmetic Fuchsian groups (see also \cite{Sarnak}).  On the other hand, Zograf \cite{Zograf} proved that 
\[ \lambda_1(\Gamma) < \frac{4(g+1)}{A} \]
if $A \geq 16(g+1)$.  From these, it follows that there are only finitely many congruence arithmetic Fuchsian groups of bounded genus, since there are only finitely many such groups of fixed signature \cite{Takeuchi} and hence bounded area.

In particular, putting these together, we obtain the following bound.

\begin{lem}[Selberg-Zograf bound] \label{SZbound}
If $\Gamma$ is a congruence arithmetic Fuchsian group of genus $g$ and area $A$, then
\[ A < \frac{64}{3}(g+1). \]
\end{lem}

From this lemma and (\ref{shimizu}) it then follows that
\begin{equation} \label{areabound}
\frac{d_F^{3/2}}{(2\pi)^{2n}}\zeta_F(2)\Phi(\frakD)\Psi(\frakN) < \frac{16}{3}(g+1),
\end{equation}
and from the trivial bound $\zeta_F(2)\Phi(\frakD)\Psi(\frakN) \geq 1$ we conclude
\begin{equation} \label{SZ}
\delta_F < (2\pi)^{4/3} \left(\frac{16(g+1)}{3}\right)^{2/(3n)} < 11.595 \cdot 6.350^{1/n}
\end{equation}
for $g \leq 2$.

Now, by the unconditional Odlyzko bounds \cite{Odlyzko}, as calculated by Martinet \cite{Martinet} (see also Cohen-Diaz y Diaz-Olivier \cite{CDO}), we have for $n \geq 11$ and $F$ totally real that $\delta_F > 14.083$, but by (\ref{SZ}) we have $\delta_F<13.717$, a contradiction.  The bounds for the remaining degrees $n \leq 10$ are summarized in the following table.

\begin{table}[h]
\begin{center}
\begin{tabular}{c|ccc|c}
$n$ & \textsf{Selberg-Zograf} & \textsf{Odlyzko} & \textsf{Odlyzko (GRH)} & $\Delta$ \\
\hline
2\rule{0pt}{2.5ex} & $<$ 29.216 & $>$ 2.223 & ($>$ 2.227) & 30 \\
3 & 21.470 & 3.610 & (3.633) & 25 \\ 
4 & 18.405 & 5.067 & (5.127) & 20 \\
5 & 16.780 & 6.523 & (6.644) & 17 \\
6 & 15.778 & 7.941 & (8.148) & 16 \\
7 & 15.098 & 9.301 & (9.617) & 15.5 \\ 
8 & 14.608 & 10.596 & (11.042) & 15 \\
9 & 14.238 & 11.823 & (12.418) & 14.5 \\
10 & 13.949 & 12.985 & (13.736) & 14 \\
\end{tabular} \\
\ \\
\textbf{Table 1.2}: Degree and root discriminant bounds
\end{center}
\end{table}

In Table 1.2, for each degree $2 \leq n \leq 10$, we list the Selberg-Zograf bound (\ref{SZ}), the unconditional Odlyzko bound, and the GRH-conditional Odlyzko bound (for comparison only).  

Let $NF(n)$ denote the set of (totally real) fields $F$ of degree $n$ up to isomorphism that satisfy the Selberg-Zograf bound with $g=2$.  All totally real fields $F$ of degree $n$ with $\delta_F \leq \Delta(n)$ have been enumerated \cite{Voight}.  From this list of fields, we find that $\#\bigcup_n NF(n)=3409$.  In Table 1.3, we list for each degree $n \leq 10$ the number $\#NF(n)$ and the smallest and largest discriminant $d_F$ for $F \in NF(n)$.  Note that $NF(10)=\emptyset$; indeed, it is conjectured \cite{Voight} that the minimal totally real field of degree $10$ has root discriminant $\delta_F \approx 14.613$, well above the Selberg-Zograf bound.

\begin{table}[h] \label{table:data}
\begin{center}
\begin{tabular}{c|c|cc}
$n=[F:\Q]$ &\ $\#NF(n)$\ \ &\ Minimal $d_F$\ &\ Maximal $d_F$ \\
\hline
\rule{0pt}{2.5ex} 2 \rule{0pt}{2.5ex} & 257 & 5 & 849 \\
3 & 377 & 49 & 9869 \\
4 & 1052 & 725 & 114629 \\
5 & 624 & 14641 & 1326492 \\
6 & 736 & 300125 & 15381949 \\
7 & 225 & 20134393 & 178799897 \\
8 & 124 & 282300416 & 2069158249 \\
9 & 14 & 9685993193 & 23041013113 \\
10 & 0 & - & - \\
\hline
\rule{0pt}{2.5ex} 
Total \rule{0pt}{2.5ex} & 3409
\end{tabular} \\
\ \\
\textbf{Table 1.3}: Totally real fields $F$ that satisfy the Selberg-Zograf bound
\end{center}
\end{table}

\section{Embedding numbers}

From the Riemann-Hurwitz formula (\ref{RH}), we see that the genus $g$ is determined by the area $A$ and the number of elliptic cycles $e_q$ (and parabolic cycles $e_\infty$ in the case of elliptic modular curves).  The numbers $e_q$, in turn, are determined by embedding numbers of commutative orders into the quaternion order $\calO$.  In this section, we bound the number of elliptic cycles $e_q$ for $\Gamma_0^\frakD(\frakN)$ in terms of class numbers.  The results stated within can be found in Vign\'eras \cite[pp.~94--98]{Vigneras} and Schneider \cite[\S 3]{Schneider}.  We continue the notation introduced in \S 1; we abbreviate $\calO=\calO_0(\frakN)$ for an Eichler order of level $\frakN$.

Let $q \in \Z_{\geq 2}$.  Suppose that $e_q>0$, so that $\calO_1^*$ has a finite subgroup of order $2q$.  Then, in particular, the field $K_q=F(\zeta_{2q})$ embeds in the quaternion algebra $B$, where $\zeta_{2q}$ is a primitive $2q$th root of unity.  Such an embedding $K_q \hookrightarrow B$ exists if and only if $[K_q:F \cap K_q]=2$ and no prime $\frakp \mid \frakD=\disc(B)$ splits in $K_q$.  Moreover, if $\la \gamma \ra \subset \calO_1^*$ is a finite subgroup of order $2q$, then in fact we have an embedding $R=F(\gamma) \cap \calO \hookrightarrow \calO$; such an embedding $R \hookrightarrow \calO$ with $RF \cap \calO=R$ is said to be an \emph{optimal} embedding.  Conversely, to every optimal embedding $\iota:R \hookrightarrow \calO$, where $R$ is a quadratic $\Z_F$-order having exactly $2q$ roots of unity $w(R)=2q$, we associate the finite subgroup $\iota(R)_{\text{tors}}^* \subset \calO_1^*$.  Thus there is a bijection 
\begin{center}
$\{\text{Elliptic cycles of $\Gamma$ of order $2q$}\}$ \\
$\updownarrow$ \\
$\{\text{$\calO_1^*$-conjugacy classes of optimal embeddings $\iota:R \hookrightarrow \calO$ with $w(R)=2q$}\}$.
\end{center}

For a quadratic $\Z_F$-order $R$, we define 
\[ m(R,\calO)=\#\{\text{$\calO_1^*$-conjugacy classes of optimal embeddings $\iota:R \hookrightarrow \calO$}\}. \]  
We study the number $m(R,\calO)$ by considering the corresponding collection of local embeddings.  For a prime $\frakp$ of $\Z_F$, denote by $R_\frakp$ and $\calO_\frakp$ the localizations of $R$ and $\calO$ at $\frakp$, and let
\[ m(R_\frakp,\calO_\frakp)=\#\{\text{$\calO_\frakp^*$-conjugacy classes of optimal embeddings $\iota:R_\frakp \hookrightarrow \calO_\frakp$}\}. \]
Let $h(R),h(F)$ denote the class number of $R$ and $\Z_F$, respectively.  We then have the following formula.

\begin{lem} \label{embedclassno}
We have
\[
e_q=\frac{1}{2 h(F)}\sum_{\substack{R \subset K_q \\ w(R)=2q}} \frac{h(R)}{Q(R)} \prod_\frakp m(R_\frakp,\calO_\frakp) \]
where $Q(R)=[\N_{K_q/\Q}(R^*):\Z_F^{*2}]$.
\end{lem}

The product in this lemma makes sense because almost all factors are $1$ by Proposition \ref{embednum} below.

\begin{cor}
The signature of a Shimura curve $X_0^{\frakD}(\frakN)$ depends only on the discriminant $\frakD$ and level $\frakN$.
\end{cor}

\begin{proof}
The curve $X_0^{\frakD}(\frakN)$ up to isomorphism depends only on the group $\Gamma=\Gamma_0^{\frakD}(\frakN)$ up to conjugation in $PSL_2(\R)$, which depends only on the quaternion algebra $B$.  Now a quaternion algebra $B$ of discriminant $\frakD$ is uniquely defined up to the choice of the unique (real) infinite split place.  Hence the corollary is true if $F=\Q$, so we may assume that $\Gamma$ has no parabolic cycles.  We then note that the formula for the number of elliptic cycles in Lemma \ref{embedclassno} is independent of the choice of real place, as is area by the formula (\ref{shimizu}), and hence also the genus.
\end{proof}

It is well-known that the quantities in Lemma \ref{embedclassno} are computable.  For example, we have
\[ h(R)=\frac{h(K_q)}{[\Z_{K_q}^*:R^*]} N(\frakf) \prod_{\frakp \mid \frakf} \left(1-\legen{K_q}{\frakp}\frac{1}{N\frakp}\right) \]
where $\frakf$ denotes the conductor of $R$ in $\Z_{K_q}$ and
\[ \legen{K}{\frakp}=
\begin{cases}
-1, & \text{if $\frakp$ is inert in $K$}; \\
0, & \text{if $\frakp$ is ramified in $K$}; \\
1, & \text{if $\frakp$ splits in $K$}.
\end{cases} \]
There are methods to compute the relative class number $h(K_q)/h(F)$ numerically due to Louboutin \cite{Louboutin}.  However, we also need to compute the unit indices $Q(R)$ for suborders $R \subset \Z_{K_q}$, and general methods for computing $\Z_{K_q}^*$ compute the class group $h(K_q)$ simultaneously \cite{Cohen2}; absent a strategy to recover these unit indices, we rely upon these standard methods.  In any case, the relative class numbers and unit indices need only be computed once for each field $F$.

We now give a formula for the number of local embeddings $m(R_\frakp,\calO_\frakp)$.  First we introduce some more notation.  Let $R_\frakp=\Z_{F,\frakp}[\gamma_\frakp]$ and let $f_\frakp(x)=x^2-t_\frakp x+n_\frakp$ denote the minimal polynomial of $\gamma_\frakp$, let $d_\frakp=t_\frakp^2-4n_\frakp$, and let $k(\frakp)$ denote the residue class field of $\frakp$.  


\begin{prop}\ \label{embednum}
\begin{enumalph}
\item If $\frakp \nmid \frakD\frakN$, then $m(R_\frakp,\calO_\frakp)=1$. 
\item If $\frakp \mid \frakD$, then
\[ m(R_\frakp,\calO_\frakp) =
\begin{cases}
0, & \text{if $\frakp \mid \frakf$}; \\
\displaystyle{1-\legen{K_q}{\frakp}}, & \text{otherwise}.
\end{cases} \]
\item Suppose $\frakp \mid \frakN$, let $e=\ord_\frakp(\frakN)$, and let
\[ E(r) = \left\{x \in R/\frakp^r : f_\frakp(x) \equiv 0 \psmod{\frakp^{e}}\right\}. \]
If $\ord_\frakp(d_\frakp)=0$ then 
\[ m(R_\frakp,\calO_\frakp)=\#E(e). \]
Otherwise, 
\[ m(R_\frakp,\calO_\frakp)=\#E(e)+\#\img\left(E(e+1) \to R/\frakp^{e}\right). \]
\end{enumalph}
\end{prop}

\begin{cor} \label{Neq1}
If $\frakN=\Z_F$ and $\frakD$ is prime to the conductor of $\Z_F[\zeta_{2q}]$, then
\[ e_q=\frac{1}{2 h(F)}\prod_{\frakp \mid \frakD} \left(1-\legen{K_q}{\frakp}\right)\sum_{\substack{R \subset K_q \\ w(R)=2q}} \frac{h(R)}{Q(R)}. \]
\end{cor}

For odd primes, we can make part (c) of the above proposition explicit.

\begin{lem}
If $\frakp$ is an odd prime, and $e=\ord_\frakp(\frakN) \geq 1$ and $f=\ord_\frakp(\frakf)$, then:
\begin{itemize}
\item If $\ord_\frakp(d_\frakp)=0$, then 
\[ m(R_\frakp,\calO_\frakp) = 
\begin{cases}
1, & \text{ if $e \leq f$;} \\
2, & \text{ if $f < e \leq 2f$;} \\
1+\displaystyle{\legen{d_\frakp}{\frakp}}, & \text{ if $e>2f$.} \\
\end{cases} \]
\item If $e < \ord_\frakp(d_\frakp)$, then
\[ m(R_\frakp,\calO_\frakp) = 
\begin{cases}
2, & \text{ if $e=1$;} \\
1+\#k(\frakp), & \text{ if $e=2$;} \\
e, & \text{ if $e\geq 3$ is odd;} \\
e-1, & \text{ if $e \geq 4$ is even}. \\
\end{cases} \]
\item If $e = \ord_\frakp(d_\frakp) > 0$, then
\[ m(R_\frakp, \calO_\frakp) =
\begin{cases}
\max(1,\#k(\frakp)\lfloor e/2 \rfloor), & \text{ if $e$ is odd}; \\
1+\displaystyle{\legen{K_q}{\frakp}}, & \text{ if $e$ is even}.
\end{cases} \]
\item If $e > \ord_\frakp(d_\frakp) > 0$, then
\[ m(R_\frakp, \calO_\frakp) = 
\begin{cases}
0, & \text{ if $e$ is odd}; \\
1+\displaystyle{\legen{K_q}{\frakp}}, & \text{ if $e$ is even}.
\end{cases} \]
\end{itemize}
\end{lem}

\begin{proof}
The lemma follows from a case-by-case analysis of Proposition \ref{embednum}(c), and the details are left to the reader.
\end{proof}

For even primes $\frakp$, the calculation of the local embedding number is quite subtle and as such is not conducive to a compact formula; nevertheless, the representation in Proposition \ref{embednum}(c) shows that $m(R_\frakp,\calO_\frakp)$ is effectively computable.  

\begin{rmk}
For even primes $\frakp$, we can improve upon the obvious method of calculating $m(R_\frakp,\calO_\frakp)$ by calculating $E(r)$ as follows.  The map $R/2R \to R/2R$ given by $x \mapsto x^2-t_\frakp x$ is $\F_2$-linear, and therefore by linear algebra over $\F_2$ one can compute all solutions to $f(x)=x^2-t_\frakp x+n_\frakp \equiv 0 \pmod{\frakp^r}$ if $r \leq \ord_\frakp(2)$.  For each of these solutions, one can then use a suitable Hensel lifting procedure to test which give rise to solutions modulo $\frakp^r$ for $r \leq \ord_\frakp(4)$; and then Hensel's lemma implies that each such solution lifts to a unique solution modulo $\frakp^r$ whenever $r > \ord_\frakp(4)$.  
\end{rmk}

\section{Enumerating Shimura curves}

In this section, we show how to enumerate the set $SC(F,g)$ of all Shimura curves of genus at most $g \in \Z_{>0}$ over a totally real field $F$.  For simplicity, we exclude throughout the subset of elliptic modular curves, having $F=\Q$ and $\frakD=(1)$: methods for their enumeration are well-known.  

We define the \emph{primitive area} by
\[ \Aprim = \frac{4}{(2\pi)^{2n}}d_F^{3/2} \zeta_F(2). \]
From Riemann-Hurwitz (\ref{RH}), we see that $\Aprim \in \Q$ has denominator bounded by the least common multiple of all $q$ such that $[F(\zeta_{2q}):F]=2$ (which in particular requires that $F$ contains the totally real subfield $\Q(\zeta_{2q})^+$ of $\Q(\zeta_{2q})$).  Therefore, $\Aprim$ can be computed from the usual Dirichlet series or Euler product expansion for $\zeta_F(2)$ with the required precision; see Dokchitser \cite{Dokchitser}.

Already from the bound (\ref{areabound}), we obtain
\begin{equation} \label{boundgmax}
\Phi(\frakD)\Psi(\frakN) < \frac{64}{3\Aprim}(g+1) = M(F,g) \in \Q.
\end{equation}
Hence to enumerate the set $SC(F,g)$, it suffices to list the finitely many possibilities for $\frakD$ and then $\frakN$ satisfying the bound (\ref{boundgmax}), and for each pair to test if $g(X_0^{\frakD}(\frakN)) \leq g$.  

Although effective, this na\"ive approach is needlessly inefficient, and we improve on it as follows.  Let $\frakD \subset \Z_F$ be a discriminant.  We define $A(\frakD)=\area(X_0^{\frakD}(1))$, which by formula (\ref{shimizu}) is given by 
\begin{equation} \label{duhAD}
A(\frakD)=\Aprim\Phi(\frakD).
\end{equation}
Similarly, we define $e_q(\frakD)$ to be the number of elliptic cycles of order $q$ in $\Gamma_0^{\frakD}(1)$, which by Corollary \ref{Neq1} is given by
\begin{equation} \label{duheqD}
e_q(\frakD)=\frac{1}{2h(F)} \prod_{\frakp \mid \frakD} \left(1-\legen{K_q}{\frakp}\right) \sum_{\substack{R \subset K_q \\ w(R)=2q \\ \gcd(\frakf(R),\frakD)=1}} \frac{h(R)}{Q(R)}.
\end{equation}
Finally, let $g(\frakD)=g(X_0^{\frakD}(1))$ denote the genus of $X_0^{\frakD}(1)$.

Now suppose that $\frakd$ is a discriminant with $\frakd \mid \frakD$.  By (\ref{duhAD}), we see that $A(\frakD)=A(\frakd)\Phi(\frakD/\frakd)$.  Similarly, by (\ref{duheqD}) we have 
\begin{equation} \label{taueq0}
e_q(\frakD) \leq e_q(\frakd) \prod_{\frakp \mid (\frakD/\frakd)} \left(1-\legen{K_q}{\frakp}\right)
\end{equation}
and equality holds in (\ref{taueq0}) if and only if $e_q(\frakd)=0$ or $\gcd(\frakf(\Z_F[\zeta_{2q}]),\frakD/\frakd)=1$.  We conclude that
\begin{equation} \label{taueq}
e_q(\frakD) \leq 2^{\tau(\frakD/\frakd)} e_q(\frakd)
\end{equation}
where $\tau(\fraka)$ denotes the number of prime divisors of an ideal $\fraka \subset \Z_F$; equality holds in (\ref{taueq}) if additionally $\displaystyle{\legen{K_q}{\frakp}}=-1$ for all $\frakp \mid (\frakD/\frakd)$.  

By Riemann-Hurwitz, we have
\begin{equation} \label{gbound}
\begin{array}{rl}
2g(\frakD)-2 &= A(\frakD)-\sum_q e_q(\frakD)(1-1/q) \\
&\geq A(\frakd)\Phi(\frakD/\frakd) - 2^{\tau(\frakD/\frakd)} \sum_q e_q(\frakd)(1-1/q).  
\end{array}
\end{equation}
It follows that if $g(\frakD) \leq g$, since $g \geq 1$ we have 
\[ 2^{\tau(\frakD/\frakd)}(2g-2) \geq 2g-2 \geq 2g(\frakD)-2 \geq A(\frakd)\Phi(\frakD/\frakd)+2^{\tau(\frakD/\frakd)}\left(2g(\frakd)-2-A(\frakd)\right) \]
which after dividing by $2^{\tau(\frakD/\frakd)}A(\frakd)$ and rearranging terms becomes
\begin{equation} \label{phibound}
\frac{\Phi(\frakD/\frakd)}{2^{\tau(\frakD/\frakd)}}=\prod_{\frakp \mid (\frakD/\frakd)}\frac{\N\frakp-1}{2} \leq 1+2\left(\frac{g-g(\frakd)}{A(\frakd)}\right).
\end{equation}

For any prime $\frakp$ with $N(\frakp)>2$ we have trivially $\Phi(\frakp) \geq 2$.  Hence if $\frakD/\frakd$ is a product of primes with norm $>2$ we have
from (\ref{gbound}) that 
\[ 2g(\frakD)-2 \geq 2^{\tau(\frakD/\frakd)}(2g(\frakd)-2) \] 
and hence $g(\frakD) \geq g(\frakd)$.  

\begin{rmk}
This statement is sharp in the sense that one cannot remove the norm two condition.  For example, consider the unique cubic field $F \subset \Q(\zeta_{31})$, i.e., the totally real cubic abelian field of discriminant $961=31^2$.  Then $2$ splits completely in $F$.  We compute that $X_0^{(1)}(1)$ has signature $(2;2^4,3)$, but if $\frakD$ is the product of two primes of norm $2$, then $X_0^{\frakD}(1)$ has signature $(1;2^4,3^4)$.
\end{rmk}

Now let 
\[
\calD=\begin{cases}
\{\Z_F\}, & \text{if $d=[F:\Q]$ is odd}; \\
\{\frakp: \Phi(\frakp) \leq M(F,g)\}, & \text{if $d$ is even}.
\end{cases} \]

We begin then by computing the set $\calD$, which is effectively computable.  We initialize $S=\emptyset$.  For each $\frakd \in \calD$, we compute $A(\frakd)$ and $g(\frakd)$; note that the class numbers and local embedding numbers $m(R_\frakp,\calO_\frakp)$ need only be computed once for each field $F$.  For those $\frakd$ with $g(\frakd) \leq g$, we add $\frakd$ to $S$.  Then, for each $\frakd \in \calD$, we find all discriminants $\frakD$ with $\frakd \mid \frakD$ that satisfy the bound (\ref{phibound}).  Note in particular that if $g(\frakd) > g$, then $\frakD/\frakd$ must be a product of prime ideals of norm $2$.  For each such $\frakD$, if $g(\frakD) \leq g$ we add $\frakD$ to $S$.  Note that $X=X_0^{\frakD}(1)$ has $X \in SC(F,g)$ if and only if $\frakD \in S$.  Finally, for each $\frakD \in S$, we find all possible levels $\frakN$ satisfying (\ref{duhAD}), ordering this computation by divisibility since $g(X_0^{\frakD}(\frakM)) \leq g(X_0^{\frakD}(\frakN))$ whenever $\frakM \mid \frakN$.  This concludes the computation of $SC(F,g)$.  

\begin{rmk}
One can further improve on this method when the bound (\ref{duhAD}) is large as follows.  The set $\calD$ consists of discriminants of the form $\frakD=\frakp \fraka$ where $\fraka$ is a squarefree product of ideals of norm $2$ and $\frakp$ is either a prime ideal or equal to $\Z_F$.  There are as such only finitely many possibilities for $\fraka$.  For each such we characterize the possibilities for $\frakp$ as follows: by (\ref{taueq0}), if $\frakp \nmid q\Z_F$ then 
\[ e_q(\frakD)=\left(1-\legen{K_q}{\frakp}\right)e_q(\fraka). \]
If $\frakp$ is such that $\sigma_q=1-(K_q/\frakp) \in \{0,2\}$ then 
\begin{equation} \label{ADphip}
A(\frakD)=\Aprim\Phi(\frakp)=2g(\frakD)-2+\sum_q \sigma_q e_q(\fraka) \left(1-\frac{1}{q}\right).
\end{equation}
Let $Q=\{q: e_q(\fraka)>0\}$.  Then for each $(\sigma_q) \in \{0,2\}^Q$, we can substitute into equation (\ref{ADphip}) with $g(\frakD) \leq g$ and then solve for $\Phi(\frakp)$ (and hence $\frakp$).  In particular, the right-hand side of (\ref{ADphip}) must be an integer of the form $N\frakp + 1$, and the corresponding prime $\frakp$ must have the specified values of the Artin symbol.  This idea cuts down dramatically on the number of discriminants that must be tested to compute $\calD$ when the bound (\ref{duhAD}) is large.

We also note that if $X_0^{\frakD}(\frakN)$ has genus at least 2, then for any multiple $\frakN'$ of $\frakN$ the curve $X_0^{\frakD}(\frakN')$ covers $X_0^{\frakD}(\frakN)$ and therefore has genus at least $3$.
\end{rmk}

\section{Tables of Shimura curves of genus at most two}



In the following tables, for each Shimura curve $X_0^{\frakD}(\frakN)$ of genus $g \leq 2$, we record the discriminant $d_F$ of $F$, the norms $D=\N_{F/Q}(\frakD)$, $N=\N_{F/\Q}(\frakN)$, and the signature $\sigma$ of $\Gamma_0^{\frakD}(\frakN)$.  This way of recording curves is compact but ambiguous; nevertheless, in all but a handful of cases, the field $F$ is determined by its discriminant, and for \emph{any} choice of squarefree $\frakD$ (including the choice of ramified infinite places) and coprime $\frakN$, the curve $X_0^{\frakD}(\frakN)$ has the given signature.  For the handful of exceptions, we refer to the complete tables that are available online \cite{Voightonline}.

\renewcommand{\baselinestretch}{0.95}

\begin{center}
\begin{footnotesize}
\[
\begin{array}{ccc|ccc|ccc}
D & N & \sigma & D & N & \sigma & D & N & \sigma \\
\hline
1 & 1 & \T(0; 2,3; 1) &	1 & 19 & (1; 3^2; 2) &	6 & 1 & (0;2^2,3^2) \\
 & 2 & (0; 2; 2) &	 & 20 & (1; -; 6) &	  & 5 & (1;2^4) \\
 & 3 & (0; 3; 2) &	 & 21 & (1; 3^2; 4) &	  & 7 & (1;3^4) \\
 & 4 & (0; -; 3) &	 & 22 & (2; -; 4) &	  & 13 & (1;2^4,3^4) \\
 & 5 & (0; 2^2; 2) &	 & 23 & (2; -; 2) &	10 & 1 & (0;3^4) \\
 & 6 & (0; -; 4) &	 & 24 & (1; -; 8) &	   & 3 & (1;3^4) \\
 & 7 & (0; 3^2; 2) &	 & 25 & (0; 2^2; 6) &	   & 7 & (1;3^8) \\
 & 8 & (0; -; 4) &	 & 26 & (2; 2^2; 4) &	14 & 1 & (1;2^2) \\
 & 9 & (0; -; 4) &	 & 27 & (1; -; 6) &	15 & 1 & (1;3^2) \\
 & 10 & (0; 2^2; 4) &	 & 28 & (2; -; 6) &	21 & 1 & (1;2^4) \\
 & 11 & (1; -; 2) &	 & 29 & (2; 2^2; 2) &	22 & 1 & (0;2^2,3^4) \\
 & 12 & (0; -; 6) &	 & 31 & (2; 3^2; 2) &	26 & 1 & (2;-) \\
 & 13 & (0; 2^2, 3^2; 2) &	 & 32 & (1; -; 8) &	33 & 1 & (1;2^4,3^2) \\
 & 14 & (1; -; 4) &	 & 36 & (1; -; 12) &	34 & 1 & (1;3^4) \\
 & 15 & (1; -; 4) &	 & 37 & (2; 2^2,3^2; 2) &	38 & 1 & (2;2^2) \\
 & 16 & (0; -; 6) &	 & 49 & (1; 3^2; 8) &	46 & 1 & (1;2^2,3^4) \\
 & 17 & (1; 2^2; 2) & 	 & 50 & (2; 2^2; 12) &	58 & 1 & (2;3^4) \\
 & 18 & (0; -; 8) &		
 \end{array}
\] \\
\end{footnotesize}
\textbf{Table 4.1}: Shimura curves with $F=\Q$
\end{center}

\newpage

\begin{center}
\begin{footnotesize}
\[
\begin{array}{cccc|cccc|cccc|cccc}
d_F\hspace{-1ex} & D & N & \sigma & d_F & D & N & \sigma & d_F & D & N & \sigma & d_F & D & N & \sigma \\
\hline
5 & 4 & 1 & \T(0;2,5^2) &	8 & 2 & 1 & (0;3^2,4) &	17 & 2 & 1 & (0;2^2,3^2) &	40 & 2 & 1 & (0;2,3^4) \\
 & 4 & 5 & (0;2^2,5^2) &	 & 2 & 7 & (0;3^4) &	 & 2 & 2 & (1;2^2) &	 & 2 & 3 & (2;3^4) \\
 & 4 & 9 & (1;2^2) &	 & 2 & 9 & (0;3^2,4^2) &	 & 2 & 4 & (2;-) &	 & 3 & 1 & (0;2^6,3^2) \\
 & 4 & 11 & (0;5^4) &	 & 2 & 17 & (1;4^2) &	 & 2 & 9 & (1;2^4,3^2) &	 & 3 & 2 & (2;2^{10}) \\
 & 4 & 19 & (2;-) &	 & 2 & 23 & (2;-) &	 & 2 & 13 & (1;2^4,3^4) &	 & 5 & 1 & (2;3^4) \\
 & 4 & 25 & (2;2^2) &	 & 2 & 25 & (0;3^4,4^2) &	 & 9 & 1 & (2;3) &	 & 18 & 1 & (1;2^4,3^4) \\
 & 4 & 29 & (2;2^2) &	 & 2 & 31 & (1;3^4) &	 & 36 & 1 & (1;3^4) &	41 & 2 & 1 & (0;2^4,3^2) \\
 & 4 & 31 & (1;5^4) &	 & 2 & 41 & (2;4^2) &	 & 68 & 1 & (1;3^8) &	 & 2 & 2 & (2;2^4) \\
 & 4 & 41 & (1;2^2,5^4) &	 & 2 & 49 & (1;3^8) &	21 & 3 & 1 & (0;2^4,3) &	 & 20 & 1 & (1;3^8) \\
 & 4 & 61 & (2;2^2,5^4) &	 & 2 & 49 & (2;3^4) &	 & 3 & 4 & (1;2^4,3^2) &	44 & 2 & 1 & (0;2,3^4) \\
 & 5 & 1 & (0;3^2,5) &	 & 2 & 73 & (2;3^4,4^2) &	 & 3 & 5 & (1;2^8) &	 & 5 & 1 & (2;3^4) \\
 & 5 & 4 & (0;3^4) &	 & 7 & 1 & (0;2^2,4^2) &	 & 4 & 1 & (1;2^2) &	 & 7 & 1 & (2;2^{10}) \\
 & 5 & 9 & (1;3^2) &	 & 7 & 2 & (0;2^4,4^2) &	 & 5 & 1 & (0;3^5) &	53 & 4 & 1 & (2;2^3) \\
 & 5 & 11 & (1;5^2) &	 & 7 & 4 & (1;2^6) &	 & 5 & 3 & (1;3^8) &	 & 11 & 1 & (2;2^6,3^{10}) \\
 & 5 & 16 & (1;3^4) &	 & 7 & 9 & (1;2^4,4^4) &	 & 5 & 4 & (1;3^{10}) &	56 & 2 & 1 & (0;2^2,3^4) \\
 & 5 & 19 & (1;3^4) &	 & 9 & 1 & (1;3) &	 & 7 & 1 & (1;2^4) &	57 & 2 & 1 & (0;2^2,3^5) \\
 & 5 & 31 & (1;3^4,5^2) &	 & 9 & 2 & (2;-) &	 & 17 & 1 & (2;3^5) &	 & 3 & 1 & (2;2^4,3) \\
 & 9 & 1 & (0;3,5^2) &	 & 17 & 1 & (1;3^2) &	24 & 2 & 1 & (0;2,3^3) &	 & 12 & 1 & (1;2^4,3^4) \\
 & 9 & 4 & (1;3^2) &	 & 23 & 1 & (0;2^2,3^2,4^2) &	 & 2 & 3 & (0;3^6) &	60 & 2 & 1 & (0;3^6) \\
 & 9 & 5 & (1;5^2) &	 & 23 & 2 & (2;2^4,4^2) &	 & 2 & 5 & (2;2^2) &	 & 2 & 3 & (1;3^{12}) \\
 & 9 & 11 & (1;5^4) &	 & 25 & 1 & (2;-) &	 & 2 & 9 & (0;3^{12}) &	 & 3 & 1 & (1;2^8) \\
 & 11 & 1 & (0;2^2,3^2) &	 & 31 & 1 & (1;2^2,4^2) &	 & 3 & 1 & (0;2^6) &	61 & 3 & 1 & (0;2^6,3^4) \\
 & 11 & 4 & (0;2^2,3^4) &	 & 41 & 1 & (2;3^2) &	 & 3 & 2 & (0;2^{10}) &	 & 5 & 1 & (2;3^8) \\
 & 11 & 5 & (1;2^4) &	 & 47 & 1 & (1;2^2,3^2,4^2) &	 & 3 & 4 & (1;2^{12}) &	65 & 2 & 1 & (0;2^4,3^4) \\
 & 11 & 9 & (1;2^4,3^2) &	 & 71 & 1 & (2;2^2,3^2,4^2) &	 & 3 & 5 & (1;2^{12}) &	 & 20 & 1 & (1;3^{16}) \\
 & 19 & 1 & (0;2^2,5^2) &	 & 98 & 1 & (1;4^4) &	 & 5 & 1 & (1;3^3) &	69 & 3 & 1 & (1;2^8) \\
 & 19 & 4 & (2;2^2) &	12 & 2 & 1 & (0;3^2,6) &	 & 50 & 1 & (1;3^{12}) &	 & 5 & 1 & (2;3^9) \\
 & 19 & 5 & (1;2^4,5^2) &	 & 2 & 3 & (0;3^4) &	28 & 2 & 1 & (0;3^4) &	73 & 2 & 1 & (1;2^2,3^4) \\
 & 29 & 1 & (0;3^2,5^2) &	 & 2 & 9 & (0;3^6) &	 & 2 & 3 & (1;3^4) &	 & 12 & 1 & (1;2^4,3^8) \\
 & 29 & 4 & (2;3^4) &	 & 2 & 11 & (2;-) &	 & 2 & 7 & (1;3^8) &	76 & 2 & 1 & (1;2,3^4) \\
 & 31 & 1 & (1;2^2) &	 & 2 & 13 & (0;3^4,6^2) &	 & 3 & 1 & (0;2^4,3^2) &	 & 3 & 1 & (1;2^{10},3^2) \\
 & 41 & 1 & (1;3^2) &	 & 2 & 25 & (1;3^4,6^2) &	 & 3 & 2 & (1;2^8) &	85 & 3 & 1 & (1;2^4,3^6) \\
 & 49 & 1 & (1;5^2) &	 & 2 & 37 & (2;3^4,6^2) &	 & 3 & 4 & (2;2^{12}) &	88 & 2 & 1 & (0;2,3^8) \\
 & 59 & 1 & (0;2^2,3^2,5^2) &	 & 3 & 1 & (0;2^3,6) &	 & 7 & 1 & (2;2^4) &	 & 3 & 1 & (2;2^6,3^4) \\
 & 61 & 1 & (2;-) &	 & 3 & 2 & (0;2^6) &	 & 18 & 1 & (1;3^4) &	89 & 2 & 1 & (1;2^6,3^2) \\
 & 71 & 1 & (1;2^2,3^2) &	 & 3 & 4 & (0;2^8) &	29 & 4 & 1 & (1;2^3) &	92 & 2 & 1 & (0;3^8) \\
 & 79 & 1 & (1;2^2,5^2) &	 & 3 & 8 & (1;2^8) &	 & 5 & 1 & (0;3^6) &	93 & 3 & 1 & (2;2^4,3^3) \\
 & 89 & 1 & (1;3^2,5^2) &	 & 3 & 13 & (1;2^6,6^2) &	 & 5 & 4 & (2;3^{12}) &	97 & 2 & 1 & (2;2^2,3^4) \\
 & 101 & 1 & (2;3^2) &	 & 11 & 1 & (0;2^3,3^2,6) &	 & 7 & 1 & (1;2^6) &	104 & 2 & 1 & (1;2^3,3^4) \\
 & 109 & 1 & (2;5^2) &	 & 11 & 2 & (2;2^6) &	 & 9 & 1 & (2;3^3) &	105 & 2 & 1 & (1;2^4,3^6) \\
 & 131 & 1 & (2;2^2,3^2) &	 & 13 & 1 & (2;-) &	33 & 2 & 1 & (0;2^2,3^3) &	109 & 3 & 1 & (2;2^6,3^6) \\
 & 139 & 1 & (2;2^2,5^2) &	 & 23 & 1 & (1;2^3,3^2,6) &	 & 2 & 2 & (2;2^2) &	113 & 2 & 1 & (1;2^4,3^6) \\
 & 149 & 1 & (2;3^2,5^2) &	 & 66 & 1 & (1;6^4) &	 & 2 & 3 & (1;3^6) &	120 & 2 & 1 & (0;2^2,3^{10}) \\
 & 179 & 1 & (2;2^2,3^2,5^2) &	13 & 3 & 1 & (0;2^2,3^2) &	 & 3 & 1 & (1;2^4) &	129 & 2 & 1 & (2;2^6,3^5) \\
 & 180 & 1 & (1;5^4) &	 & 3 & 3 & (1;3^2) &	 & 12 & 1 & (1;2^4) &	137 & 2 & 1 & (2;2^4,3^6) \\
&&&&	 & 3 & 4 & (0;2^2,3^4) &	 & 44 & 1 & (1;2^4,3^{12}) &	140 & 2 & 1 & (1;2^2,3^8) \\
&&&&	 & 3 & 13 & (1;2^4,3^4) &	37 & 3 & 1 & (0;2^2,3^4) &	152 & 2 & 1 & (1;2^3,3^8) \\
&&&&	 & 4 & 1 & (1;2) &	 & 3 & 4 & (2;2^2,3^8) &	156 & 2 & 1 & (2;3^{10}) \\
&&&&	 & 4 & 3 & (2;-) &	 & 4 & 1 & (2;2) &	168 & 2 & 1 & (1;2^2,3^{12}) \\
&&&&	 & 13 & 1 & (2;-) &	 & 11 & 1 & (2;2^2,3^8) &	172 & 2 & 1 & (2;2,3^{12}) \\
&&&&	 & 17 & 1 & (1;3^4) &	& \\	
&&&&	 & 23 & 1 & (1;2^2,3^4) &	& \\	
&&&&	 & 29 & 1 & (2;3^4) &		
\end{array}
\] \\
\end{footnotesize}
\textbf{Table 4.2}: Shimura curves with $[F:\Q]=2$
\end{center}

\newpage

\begin{center}
\begin{footnotesize}
\[
\begin{array}{cccc|cccc|cccc|cccc}
d_F & D & N & \sigma & d_F & D & N & \sigma & d_F & D & N & \sigma & d_F & D & N & \sigma \\
\hline
49 & 1 & 1 & \T(0;2,3,7) &	148 & 1 & 1 & (0;2^3,3) &	321 & 1 & 1 & (0;2,3^3) &	961 & 1 & 1 & (2;2^4,3) \\
 & 1 & 7 & (0;3^2,7) &	 & 1 & 2 & (0;2^5) &	 & 1 & 3 & (0;3^6) &	 & 4 & 1 & (1;2^4,3^4) \\
 & 1 & 8 & (0;2,7^2) &	 & 1 & 4 & (0;2^6) &	 & 1 & 3 & (1;3^3) &	985 & 1 & 1 & (0;2^6,3^2) \\
 & 1 & 13 & (0;2^2,3^2) &	 & 1 & 5 & (0;2^6) &	 & 1 & 7 & (1;3^6) &	993 & 1 & 1 & (0;2^3,3^5) \\
 & 1 & 27 & (1;3) &	 & 1 & 8 & (1;2^4) &	 & 1 & 9 & (0;3^{12}) &	1016 & 1 & 1 & (1;2^6,3^2) \\
 & 1 & 29 & (0;2^2,7^2) &	 & 1 & 10 & (0;2^{10}) &	 & 9 & 1 & (1;2^4) &	 & 4 & 1 & (0;2^2,3^8) \\
 & 1 & 41 & (1;2^2) &	 & 1 & 13 & (0;2^6,3^2) &	 & 24 & 1 & (2;2^2,3^6) &	 & 6 & 1 & (1;2^{12},3^4) \\
 & 1 & 43 & (0;3^2,7^2) &	 & 1 & 17 & (1;2^6) &	361 & 1 & 1 & (0;2,3^3) &	1076 & 1 & 1 & (0;2^6,3^4) \\
 & 1 & 49 & (1;3^2) &	 & 1 & 19 & (2;3^2) &	 & 1 & 7 & (1;3^6) &	 & 6 & 1 & (1;2^4,3^8) \\
 & 1 & 56 & (1;7^2) &	 & 1 & 20 & (1;2^{12}) &	404 & 1 & 1 & (0;2^3,3^2) &	1101 & 1 & 1 & (1;2^2,3^5) \\
 & 1 & 64 & (1;7^2) &	 & 1 & 25 & (1;2^6,3^2) &	 & 1 & 2 & (1;2^5) &	 & 6 & 1 & (1;2^4,3^{10}) \\
 & 1 & 71 & (1;7^2) &	 & 1 & 25 & (2;2^6) &	 & 1 & 3 & (2;3^2) &	1129 & 1 & 1 & (1;2^2,3^4) \\
 & 1 & 83 & (2;-) &	 & 1 & 26 & (2;2^{10}) &	 & 1 & 4 & (2;2^6) &	1229 & 1 & 1 & (1;2^4,3^4) \\
 & 1 & 91 & (1;3^4) &	 & 1 & 29 & (2;2^6) &	 & 6 & 1 & (0;2^2,3^4) &	 & 6 & 1 & (1;2^8,3^8) \\
 & 1 & 97 & (1;2^2,3^2) &	 & 1 & 37 & (2;2^6,3^2) &	 & 22 & 1 & (2;2^2,3^8) &	1257 & 1 & 1 & (0;2^4,3^6) \\
 & 1 & 104 & (2;2^2) &	 & 10 & 1 & (0;3^4) &	469 & 1 & 1 & (0;2^2,3^3) &	1300 & 1 & 1 & (0;2^9,3^3) \\
 & 1 & 113 & (1;2^2,7^2) &	 & 26 & 1 & (2;-) &	 & 1 & 2 & (2;2^2) &	1304 & 1 & 1 & (2;2^6,3^2) \\
 & 1 & 125 & (2;2^2) &	 & 34 & 1 & (1;3^4) &	 & 1 & 4 & (1;2^2,3^6) &	 & 4 & 1 & (1;2^2,3^8) \\
 & 1 & 127 & (1;3^2,7^2) &	 & 38 & 1 & (2;2^2) &	 & 8 & 1 & (2;2^2) &	1345 & 1 & 1 & (0;2^5,3^5) \\
 & 1 & 139 & (2;3^2) &	 & 46 & 1 & (1;2^2,3^4) &	 & 22 & 1 & (1;2^4,3^{12}) &	1369 & 1 & 1 & (1;2^3,3^3) \\
 & 1 & 169 & (1;2^4,3^4) &	 & 54 & 1 & (2;2^2,3^2) &	473 & 1 & 1 & (0;2^3,3^2) &	1373 & 1 & 1 & (1;2^6,3^4) \\
 & 1 & 169 & (2;2^2,3^2) &	 & 58 & 1 & (2;3^4) &	 & 1 & 3 & (2;3^2) &	 & 6 & 1 & (1;2^{12},3^8) \\
 & 1 & 181 & (2;2^2,3^2) &	169 & 1 & 1 & (0;2^3,3) &	 & 1 & 5 & (2;2^6) &	1384 & 1 & 1 & (2;2^6,3^2) \\
 & 1 & 197 & (2;2^2,7^2) &	 & 1 & 5 & (0;2^6) &	564 & 1 & 1 & (0;2^3,3^3) &	 & 4 & 1 & (1;2^2,3^8) \\
 & 1 & 211 & (2;3^2,7^2) &	 & 1 & 8 & (1;2^3) &	 & 1 & 2 & (2;2^5) &	1396 & 1 & 1 & (0;2^{12},3^2) \\
 & 1 & 232 & (2;2^2,7^4) &	 & 1 & 13 & (0;2^6,3^2) &	 & 1 & 3 & (2;3^6) &	1425 & 1 & 1 & (1;2^3,3^5) \\
 & 56 & 1 & (1;2^2) &	 & 1 & 25 & (1;2^{12}) &	 & 6 & 1 & (0;2^2,3^6) &	1436 & 1 & 1 & (2;2^8,3^2) \\
 & 91 & 1 & (1;7^2) &	 & 1 & 25 & (2;2^6) &	 & 6 & 1 & (2;2^2) &	 & 4 & 1 & (2;3^8) \\
 & 104 & 1 & (2;-) &	 & 25 & 1 & (1;3^4) &	 & 9 & 1 & (1;2^{12}) &	1489 & 1 & 1 & (0;2^4,3^6) \\
 & 169 & 1 & (1;7^4) &	 & 40 & 1 & (2;3^4) &	568 & 1 & 1 & (0;2^6,3) &	1492 & 1 & 1 & (1;2^6,3^4) \\
 & 189 & 1 & (1;2^4,7^2) &	229 & 1 & 1 & (0;2^2,3^2) &	 & 1 & 2 & (1;2^{10}) &	1509 & 1 & 1 & (2;2^2,3^6) \\
 & 216 & 1 & (2;2^2,3^2) &	 & 1 & 2 & (1;2^2) &	 & 1 & 2 & (2;2^6) &	1524 & 1 & 1 & (2;2^3,3^5) \\
 & 232 & 1 & (2;3^4) &	 & 1 & 4 & (0;2^2,3^4) &	 & 4 & 1 & (0;2^2,3^4) &	1556 & 1 & 1 & (2;2^6,3^2) \\
81 & 1 & 1 & (0;2,3,9) &	 & 1 & 4 & (2;-) &	621 & 1 & 1 & (0;2^2,3^4) &	1573 & 1 & 1 & (1;2^6,3^5) \\
 & 1 & 3 & (0;3^2,9) &	 & 1 & 7 & (1;3^4) &	 & 1 & 3 & (2;3^7) &	1593 & 1 & 1 & (0;2^4,3^8) \\
 & 1 & 8 & (1;2) &	 & 1 & 13 & (1;2^4,3^4) &	 & 1 & 4 & (2;2^2,3^8) &	1620 & 1 & 1 & (1;2^9,3^4) \\
 & 1 & 9 & (0;3^4) &	 & 8 & 1 & (1;2^2) &	 & 6 & 1 & (1;2^4,3^2) &	1708 & 1 & 1 & (2;2^{10},3^3) \\
 & 1 & 17 & (1;2^2) &	 & 14 & 1 & (1;2^4) &	697 & 1 & 1 & (0;2^4,3^2) &	 & 4 & 1 & (1;2^2,3^{12}) \\
 & 1 & 19 & (0;3^2,9^2) &	 & 46 & 1 & (1;2^4,3^8) &	733 & 1 & 1 & (0;2^4,3^3) &	1765 & 1 & 1 & (1;2^{10},3^4) \\
 & 1 & 24 & (2;-) &	257 & 1 & 1 & (0;2^2,3^2) &	 & 10 & 1 & (1;3^{12}) &	1825 & 1 & 1 & (1;2^6,3^4) \\
 & 1 & 27 & (1;3^3) &	 & 1 & 3 & (1;3^2) &	756 & 1 & 1 & (0;2^3,3^4) &	1901 & 1 & 1 & (2;2^6,3^6) \\
 & 1 & 37 & (0;2^2,3^2,9^2) &	 & 1 & 5 & (1;2^4) &	 & 6 & 1 & (2;2^2,3^2) &	1929 & 1 & 1 & (1;2^2,3^{10}) \\
 & 1 & 53 & (2;2^2) &	 & 1 & 7 & (1;3^4) &	761 & 1 & 1 & (0;2^2,3^4) &	1937 & 1 & 1 & (1;2^6,3^6) \\
 & 1 & 57 & (1;3^4,9^2) &	 & 1 & 8 & (2;2^2) &	785 & 1 & 1 & (0;2^5,3^2) &	1940 & 1 & 1 & (2;2^{12},3^2) \\
 & 1 & 73 & (1;2^2,3^2,9^2) &	 & 1 & 9 & (1;2^4,3^2) &	788 & 1 & 1 & (0;2^6,3^2) &	1944 & 4 & 1 & (1;2^2,3^{16}) \\
 & 1 & 109 & (2;2^2,3^2,9^2) &	 & 15 & 1 & (1;3^4) &	 & 1 & 2 & (2;2^{10}) &	1957 & 1 & 1 & (2;2^8,3^4) \\
 & 24 & 1 & (0;2^2,9^2) &	 & 21 & 1 & (1;2^8) &	 & 6 & 1 & (1;2^4,3^4) &	2057 & 1 & 1 & (2;2^6,3^4) \\
 & 51 & 1 & (1;9^2) &	 & 24 & 1 & (1;2^4,3^4) &	837 & 1 & 1 & (0;2^4,3^4) &	2177 & 1 & 1 & (2;2^5,3^6) \\
 & 57 & 1 & (1;2^4) &	316 & 1 & 1 & (0;2^4,3) &	 & 6 & 1 & (1;2^8,3^2) &	2233 & 1 & 1 & (2;2^5,3^5) \\
 & 136 & 1 & (1;3^4,9^4) &	 & 1 & 2 & (0;2^8) &	 & 10 & 1 & (1;3^{16}) &	2241 & 1 & 1 & (2;2^9,3^4) \\
&&&&	 & 1 & 2 & (1;2^4) &	892 & 1 & 1 & (0;2^8,3^2) &	& \\
&&&&	 & 1 & 4 & (0;2^{12}) &	 & 1 & 2 & (2;2^{16}) &	& \\
&&&&	 & 1 & 4 & (2;2^8) &	 & 4 & 1 & (0;3^8) &	& \\
&&&&	 & 1 & 8 & (1;2^{16}) &	940 & 1 & 1 & (0;2^{10},3) &	& \\
&&&&	 & 4 & 1 & (0;3^4) &	 & 1 & 2 & (2;2^{18}) &	& \\
&&&&	 & 22 & 1 & (1;2^8,3^4) &	 & 4 & 1 & (1;2^2,3^4) &	
\end{array}
\] \\
\end{footnotesize}
\textbf{Table 4.3}: Shimura curves with $[F:\Q]=3$
\end{center}

\newpage

\begin{center}
\begin{footnotesize}
\[
\begin{array}{cccc|cccc|cccc}
d_F & D & N & \sigma & d_F & D & N & \sigma & d_F & D & N & \sigma \\
\hline
725 & 11 & 1 & \T(0;2^2,3^2) &	2000 & 4 & 1 & (0;5^2,10) &	4352 & 2 & 1 & (0;3^4) \\
 & 11 & 16 & (2;2^2,3^4) &	 & 4 & 5 & (0;5^4,10^2) &	 & 2 & 7 & (1;3^8) \\
 & 16 & 1 & (1;2) &	 & 4 & 25 & (2;5^{14},10^2) &	 & 7 & 1 & (0;2^6,4^4) \\
 & 19 & 1 & (0;2^2,5^2) &	 & 5 & 1 & (0;3^4) &	 & 7 & 2 & (1;2^{12},4^8) \\
 & 25 & 1 & (1;5) &	 & 5 & 4 & (0;3^8) &	4400 & 4 & 1 & (0;2,5^4) \\
 & 29 & 1 & (0;3^2,5^2) &	 & 19 & 1 & (0;2^5,5^2,10) &	 & 5 & 1 & (0;3^4,5^2) \\
 & 31 & 1 & (1;2^2) &	 & 59 & 1 & (2;2^5,3^4,5^2,10) &	 & 11 & 1 & (0;2^{10},3^4) \\
 & 41 & 1 & (1;3^2) &	2048 & 2 & 1 & (0;3^2,8) &	4525 & 5 & 1 & (0;3^4,5^2) \\
 & 49 & 1 & (1;5^2) &	 & 2 & 17 & (2;8^2) &	 & 9 & 1 & (1;3^2,5^4) \\
 & 61 & 1 & (2;-) &	 & 17 & 1 & (2;3^2) &	 & 19 & 1 & (2;2^{10},5^4) \\
 & 79 & 1 & (1;2^2,5^2) &	 & 31 & 1 & (1;2^6,4^2,8^2) &	4752 & 3 & 1 & (0;2^5,6) \\
 & 81 & 1 & (2;3) &	 & 47 & 1 & (2;2^6,3^2,4^2,8^2) &	 & 3 & 4 & (1;2^{10},6^2) \\
 & 89 & 1 & (1;3^2,5^2) &	2225 & 4 & 1 & (0;2^2,5^2) &	 & 4 & 1 & (2;-) \\
 & 101 & 1 & (2;3^2) &	 & 4 & 4 & (2;2^2) &	 & 11 & 1 & (0;2^5,3^8,6) \\
 & 109 & 1 & (2;5^2) &	 & 19 & 1 & (1;2^4,5^2) &	4913 & 4 & 1 & (1;2^4) \\
 & 131 & 1 & (2;2^2,3^2) &	 & 29 & 1 & (1;3^6,5^2) &	5125 & 5 & 1 & (1;3^4) \\
 & 139 & 1 & (2;2^2,5^2) &	2304 & 2 & 1 & (0;3^2,12) &	 & 9 & 1 & (1;3^2,5^5) \\
 & 149 & 1 & (2;3^2,5^2) &	 & 2 & 9 & (0;3^4,12^2) &	 & 11 & 1 & (1;2^8,3^4) \\
 & 179 & 1 & (2;2^2,3^2,5^2) &	 & 2 & 25 & (2;3^4,12^2) &	5225 & 4 & 1 & (0;2^2,5^4) \\
1125 & 5 & 1 & (0;3^2,15) &	 & 9 & 1 & (2;-) &	 & 11 & 1 & (1;2^4,3^8) \\
 & 5 & 9 & (1;3^4) &	 & 23 & 1 & (0;2^6,3^2,4^3,12) &	5725 & 9 & 1 & (1;3^2,5^6) \\
 & 5 & 16 & (1;3^4,15^2) &	2525 & 5 & 1 & (0;3^2,5^2) &	 & 11 & 1 & (2;2^6,3^4) \\
 & 9 & 1 & (0;5^2,15) &	 & 11 & 1 & (0;2^6,3^2) &	5744 & 4 & 1 & (2;2) \\
 & 9 & 5 & (1;5^4) &	 & 16 & 1 & (2;2^3) &	 & 5 & 1 & (0;3^8) \\
 & 16 & 1 & (1;2^2) &	 & 29 & 1 & (2;3^2,5^4) &	 & 7 & 1 & (1;2^{10}) \\
 & 29 & 1 & (0;3^2,5^2,15) &	2624 & 4 & 1 & (1;4) &	6125 & 5 & 1 & (1;3^4,5) \\
 & 31 & 1 & (1;2^4) &	 & 7 & 1 & (0;2^4,4^2) &	6224 & 2 & 1 & (0;2,3^4) \\
 & 59 & 1 & (0;2^4,3^2,5^2,15) &	 & 7 & 4 & (2;2^8,4^2) &	 & 5 & 1 & (2;3^4) \\
 & 89 & 1 & (2;3^2,5^2,15) &	 & 17 & 1 & (1;3^6) &	 & 7 & 1 & (1;2^{14}) \\
1600 & 4 & 1 & (0;4,5^2) &	2777 & 2 & 1 & (0;2^2,3^2) &	6809 & 2 & 1 & (0;2^4,3^2) \\
 & 4 & 9 & (2;4^2) &	 & 2 & 8 & (2;2^2) &	7053 & 3 & 1 & (0;2^6,3^2) \\
 & 9 & 1 & (0;3^2,5^2) &	 & 8 & 1 & (1;2^2,3^2) &	7056 & 3 & 1 & (0;2^5,3^2,6) \\
 & 9 & 4 & (2;3^4) &	 & 11 & 1 & (1;2^4,3^2) &	7168 & 2 & 1 & (0;3^4,4) \\
 & 25 & 1 & (2;5) &	3600 & 4 & 1 & (0;5^4) &	 & 7 & 1 & (0;2^{12},4^6) \\
 & 31 & 1 & (1;2^4,4^2) &	 & 9 & 1 & (1;5^4) &	7225 & 4 & 1 & (1;2^4,5^2) \\
 & 41 & 1 & (2;3^4) &	 & 11 & 1 & (0;2^5,3^4,6) &	7232 & 2 & 1 & (0;2^2,3^2,4^2) \\
 & 71 & 1 & (2;2^4,3^4,4^2) &	3981 & 3 & 1 & (0;2^6) &	 & 2 & 2 & (2;2^4,4^2) \\
1957 & 3 & 1 & (0;2^2,3^2) &	 & 3 & 5 & (1;2^{12}) &	7488 & 2 & 1 & (0;3^4,6) \\
 & 3 & 7 & (1;3^4) &	 & 5 & 1 & (0;3^6) &	7537 & 2 & 1 & (0;2^2,3^4) \\
 & 3 & 16 & (2;2^2,3^4) &	 & 5 & 3 & (1;3^{12}) &	 & 3 & 1 & (1;2^4,3^2) \\
 & 7 & 1 & (1;2^2) &	 & 9 & 1 & (2;3^3) &	7600 & 4 & 1 & (1;2,5^4) \\
 & 16 & 1 & (2;2) &	4205 & 5 & 1 & (0;3^6) &	 & 11 & 1 & (2;2^{10},3^8) \\
 & 19 & 1 & (2;2^2) &	 & 7 & 1 & (1;2^6) &	& \\
 & 23 & 1 & (1;2^2,3^4) &	4225 & 4 & 1 & (0;2^4,5^2) &	& \\
 & 27 & 1 & (2;2^2,3^2) &	 & 9 & 1 & (1;3^4,5^2) &	
\end{array}
\] \\
\end{footnotesize}
\textbf{Table 4.4(a)}: Shimura curves with $[F:\Q]=4$ (Table 1 of 2, $d_F \leq 7600$)
\end{center}

\newpage

\begin{center}
\begin{footnotesize}
\[
\begin{array}{cccc|cccc}
d_F & D & N & \sigma & d_F & D & N & \sigma \\
\hline
7625 & 4 & 1 & \T(0;2^4,5^5) &	14272 & 3 & 1 & (0;2^{10},3^6) \\
 & 5 & 1 & (1;3^8) &	14336 & 2 & 1 & (0;3^8,4) \\
8000 & 4 & 1 & (0;2,5^7) &	14656 & 2 & 1 & (2;3^4) \\
 & 5 & 1 & (2;3^4,5) &	 & 3 & 1 & (1;2^{16},3^2) \\
8069 & 5 & 1 & (1;3^8) &	15188 & 2 & 1 & (1;2^6,3^4) \\
8112 & 3 & 1 & (0;2^{10}) &	 & 2 & 1 & (2;2^2,3^4) \\
8468 & 2 & 1 & (0;2^6,3^2) &	15317 & 2 & 1 & (1;2^4,3^6) \\
 & 2 & 1 & (1;2^2,3^2) &	15529 & 2 & 1 & (0;2^8,3^4) \\
 & 2 & 2 & (2;2^{10}) &	15952 & 2 & 1 & (2;2,3^4) \\
8525 & 5 & 1 & (1;3^4,5^4) &	 & 3 & 1 & (2;2^{14},3^2) \\
8789 & 5 & 1 & (2;3^6) &	16357 & 3 & 1 & (1;2^{12},3^4) \\
8957 & 3 & 1 & (0;2^6,3^4) &	16448 & 2 & 1 & (1;2^2,3^6,4^2) \\
9225 & 4 & 1 & (1;2^4,5^4) &	 & 2 & 1 & (2;2,3^4) \\
9248 & 2 & 1 & (0;2^4,3^4) &	16609 & 2 & 1 & (1;2^6,3^4) \\
 & 2 & 1 & (1;3^4) &	17069 & 3 & 1 & (2;2^{10},3^4) \\
9301 & 3 & 1 & (0;2^6,3^4) &	17417 & 2 & 1 & (0;2^{10},3^4) \\
 & 5 & 1 & (2;3^8) &	17424 & 2 & 1 & (2;2^3,3^4,6) \\
9909 & 3 & 1 & (1;2^8) &	17428 & 2 & 1 & (2;2^6,3^4) \\
 & 5 & 1 & (2;3^9) &	17609 & 2 & 1 & (1;2^6,3^4) \\
10025 & 4 & 1 & (1;2^6,5^4) &	17989 & 3 & 1 & (2;2^{10},3^4) \\
 & 5 & 1 & (1;3^{10},5^2) &	18097 & 3 & 1 & (1;2^8,3^{10}) \\
10273 & 2 & 1 & (0;2^4,3^4) &	18432 & 2 & 1 & (0;3^{11},4) \\
 & 3 & 1 & (1;2^8,3^2) &	18496 & 2 & 1 & (1;2^4,3^4,4^4) \\
10304 & 2 & 1 & (0;2^2,3^4,4^2) &	18688 & 2 & 1 & (1;3^8,4) \\
10889 & 2 & 1 & (0;2^4,3^4) &	18736 & 3 & 1 & (2;2^{10},3^6) \\
11025 & 4 & 1 & (1;2^8,5^4) &	19429 & 3 & 1 & (2;2^{10},3^6) \\
 & 5 & 1 & (2;3^9,5^2) &	19796 & 2 & 1 & (1;2^{12},3^4) \\
11197 & 3 & 1 & (1;2^4,3^4) &	20808 & 2 & 1 & (0;2^4,3^{11}) \\
11324 & 2 & 1 & (0;2^2,3^6) &	21208 & 2 & 1 & (2;2^2,3^8) \\
11344 & 2 & 1 & (1;2,3^4) &	21308 & 2 & 1 & (1;2^2,3^{10}) \\
 & 3 & 1 & (0;2^{14},3^2) &	21312 & 2 & 1 & (1;3^{12},6) \\
11348 & 2 & 1 & (0;2^6,3^4) &	21469 & 3 & 1 & (2;2^{10},3^8) \\
 & 2 & 1 & (1;2^2,3^4) &	21568 & 2 & 1 & (2;2^4,3^4,4^4) \\
11525 & 5 & 1 & (2;3^6,5^4) &	21964 & 2 & 1 & (1;3^{12}) \\
11661 & 3 & 1 & (1;2^4,3^6) &	22545 & 2 & 1 & (1;2^6,3^9) \\
12357 & 3 & 1 & (2;2^4,3^3) &	22676 & 2 & 1 & (2;2^{12},3^4) \\
12544 & 2 & 1 & (0;3^8) &	22784 & 2 & 1 & (2;3^8,4) \\
13025 & 4 & 1 & (2;2^8,5^4) &	23297 & 2 & 1 & (1;2^6,3^8) \\
13068 & 2 & 1 & (0;3^9) &	23377 & 2 & 1 & (2;2^4,3^8) \\
 & 3 & 1 & (1;2^{16}) &	23552 & 2 & 1 & (1;3^{12},4) \\
13448 & 2 & 1 & (0;2^4,3^6) &	23724 & 2 & 1 & (2;3^{11}) \\
13625 & 4 & 1 & (2;2^4,5^7) &	24417 & 2 & 1 & (1;2^6,3^9) \\
13676 & 2 & 1 & (1;3^6) &	25961 & 2 & 1 & (2;2^6,3^8) \\
13768 & 2 & 1 & (0;2^2,3^8) &	26825 & 2 & 1 & (2;2^{10},3^6) \\
 & 3 & 1 & (1;2^{12},3^4) &	26873 & 2 & 1 & (2;2^4,3^{10}) \\
13824 & 2 & 1 & (0;3^8,6) &	30056 & 2 & 1 & (1;2^4,3^{16}) \\
 & 3 & 1 & (1;2^{15},6) &	30776 & 2 & 1 & (2;2^4,3^{14}) \\
13968 & 2 & 1 & (1;2^3,3^4,6) &	& \\
14013 & 3 & 1 & (1;2^4,3^8) &	
\end{array}
\] \\
\end{footnotesize}
\textbf{Table 4.4(b)}: Shimura curves with $[F:\Q]=4$ (Table 2 of 2, $d_F > 7600$)
\end{center}

\newpage

\begin{center}
\begin{footnotesize}
\[
\begin{array}{cccc|cccc|cccc}
d_F & D & N & \sigma & d_F & D & N & \sigma & d_F & D & N & \sigma \\
\hline
14641 & 1 & 1 & \T(0;2,3,11) &	135076 & 1 & 1 & (0;2^6,3^3) &	240133 & 1 & 1 & (1;2^{10},3^3) \\
 & 1 & 11 & (1;11) &	138136 & 1 & 1 & (0;2^6,3^3) &	240881 & 1 & 1 & (0;2^7,3^6) \\
 & 1 & 23 & (1;11^2) &	138917 & 1 & 1 & (0;2^4,3^4) &	242773 & 1 & 1 & (0;2^6,3^7) \\
 & 1 & 32 & (2;2) &	144209 & 1 & 1 & (0;2^4,3^4) &	245992 & 1 & 1 & (1;2^{12},3^2) \\
 & 1 & 43 & (2;3^2) &	147109 & 1 & 1 & (0;2^8,3^2) &	246832 & 1 & 1 & (1;2^{10},3^4) \\
 & 1 & 67 & (2;3^2,11^2) &	149169 & 1 & 1 & (0;2^3,3^5) &	249689 & 1 & 1 & (0;2^{10},3^4) \\
24217 & 1 & 1 & (0;2^3,3) &	153424 & 1 & 1 & (0;2^7,3^3) &	255877 & 1 & 1 & (1;2^6,3^5) \\
 & 1 & 5 & (0;2^6) &	157457 & 1 & 1 & (0;2^5,3^4) &	265504 & 1 & 1 & (1;2^{18},3^2) \\
 & 1 & 17 & (1;2^6) &	160801 & 1 & 1 & (0;2^5,3^4) &	270017 & 1 & 1 & (2;2^5,3^4) \\
 & 1 & 25 & (2;2^6) &	161121 & 1 & 1 & (0;2^3,3^6) &	273397 & 1 & 1 & (2;2^6,3^5) \\
 & 1 & 29 & (2;2^6) &	170701 & 1 & 1 & (0;2^4,3^5) &	274129 & 1 & 1 & (1;2^9,3^3) \\
 & 1 & 37 & (2;2^6,3^2) &	173513 & 1 & 1 & (0;2^5,3^4) &	284897 & 1 & 1 & (2;2^6,3^4) \\
36497 & 1 & 1 & (0;2^2,3^2) &	176281 & 1 & 1 & (0;2^7,3^3) &	287349 & 1 & 1 & (1;2^4,3^9) \\
 & 1 & 3 & (1;3^2) &	176684 & 1 & 1 & (0;2^{10},3^2) &	288565 & 1 & 1 & (1;2^6,3^7) \\
 & 1 & 13 & (1;2^4,3^4) &	179024 & 1 & 1 & (1;2^7,3^2) &	288633 & 1 & 1 & (2;2^5,3^5) \\
38569 & 1 & 1 & (0;2^2,3^2) &	180769 & 1 & 1 & (0;2^6,3^4) &	303952 & 1 & 1 & (1;2^{14},3^4) \\
 & 1 & 7 & (1;3^4) &	181057 & 1 & 1 & (1;2^3,3^4) &	305617 & 1 & 1 & (1;2^9,3^5) \\
 & 1 & 13 & (1;2^4,3^4) &	186037 & 1 & 1 & (1;2^4,3^4) &	307145 & 1 & 1 & (2;2^8,3^4) \\
65657 & 1 & 1 & (0;2^3,3^2) &	195829 & 1 & 1 & (1;2^6,3^3) &	307829 & 1 & 1 & (1;2^8,3^6) \\
 & 1 & 3 & (2;3^2) &	202817 & 1 & 1 & (0;2^5,3^6) &	310097 & 1 & 1 & (2;2^8,3^4) \\
 & 1 & 5 & (2;2^6) &	205225 & 1 & 1 & (1;2^5,3^4) &	310257 & 1 & 1 & (2;2^3,3^9) \\
70601 & 1 & 1 & (0;2^3,3^2) &	207184 & 1 & 1 & (1;2^7,3^3) &	312617 & 1 & 1 & (2;2^5,3^6) \\
81509 & 1 & 1 & (0;2^4,3^2) &	210557 & 1 & 1 & (0;2^8,3^4) &	313905 & 1 & 1 & (2;2^7,3^5) \\
 & 1 & 2 & (2;2^4) &	216637 & 1 & 1 & (1;2^8,3^3) &	329977 & 1 & 1 & (2;2^{10},3^4) \\
81589 & 1 & 1 & (0;2^4,3^2) &	218524 & 1 & 1 & (0;2^{12},3^3) &	339509 & 1 & 1 & (2;2^6,3^8) \\
 & 1 & 2 & (2;2^4) &	220036 & 1 & 1 & (1;2^6,3^5) &	341692 & 1 & 1 & (1;2^{20},3^3) \\
89417 & 1 & 1 & (0;2^4,3^2) &	220669 & 1 & 1 & (1;2^{10},3^2) &	347317 & 1 & 1 & (2;2^8,3^6) \\
101833 & 1 & 1 & (0;2^3,3^3) &	223824 & 1 & 1 & (1;2^7,3^5) &	354969 & 1 & 1 & (2;2^7,3^8) \\
106069 & 1 & 1 & (0;2^4,3^3) &	223952 & 1 & 1 & (0;2^{14},3^2) &	356173 & 1 & 1 & (2;2^{10},3^5) \\
117688 & 1 & 1 & (0;2^6,3^2) &	224773 & 1 & 1 & (2;2^6,3^2) &	356789 & 1 & 1 & (1;2^{12},3^6) \\
 & 1 & 2 & (2;2^{10}) &	230224 & 1 & 1 & (1;2^{14},3^2) &	357977 & 1 & 1 & (2;2^{10},3^4) \\
122821 & 1 & 1 & (0;2^4,3^3) &	 & 4 & 1 & (2;2^2,3^8) &	373057 & 1 & 1 & (2;2^8,3^6) \\
124817 & 1 & 1 & (0;2^3,3^4) &	233489 & 1 & 1 & (0;2^6,3^6) &	375145 & 1 & 1 & (2;2^{11},3^5) \\
126032 & 1 & 1 & (0;2^7,3^2) &	236549 & 1 & 1 & (1;2^8,3^4) &	390625 & 1 & 1 & (2;2^5,3^{11}) \\
 & 1 & 2 & (2;2^{13}) &	& \\	
 & 6 & 1 & (2;2^2,3^4) &
\end{array}
\] \\
\end{footnotesize}
\textbf{Table 4.5}: Shimura curves with $[F:\Q]=5$
\end{center}

\newpage

\begin{center}
\begin{footnotesize}
\[
\begin{array}{cccc|cccc}
d_F & D & N & \sigma & d_F & D & N & \sigma \\
\hline
300125 & 29 & 1 & \T(2;3^2,5^2) &	1397493 & 3 & 1 & (1;2^8) \\
371293 & 13 & 1 & (2;13) &	1416125 & 5 & 1 & (2;3^6,5^2) \\
 & 27 & 1 & (2;2^6,3^2) &	1767625 & 4 & 1 & (2;2^6,5^4) \\
434581 & 13 & 1 & (2;7^2) &	1868969 & 2 & 1 & (0;2^6,3^4) \\
 & 27 & 1 & (2;2^6,3^2,7^2) &	2286997 & 3 & 1 & (1;2^{12},3^4) \\
453789 & 7 & 1 & (1;2^4) &	2323397 & 3 & 1 & (1;2^{12},3^4) \\
485125 & 9 & 1 & (1;3^2,5^2) &	2495261 & 3 & 1 & (2;2^{10},3^4) \\
 & 19 & 1 & (2;2^6,5^2) &	2501557 & 3 & 1 & (2;2^{10},3^4) \\
592661 & 7 & 1 & (1;2^6) &	2540864 & 2 & 1 & (1;2,3^8) \\
722000 & 4 & 1 & (1;2,5^2) &	2565429 & 3 & 1 & (2;2^{14},3^2) \\
810448 & 4 & 1 & (2;2) &	2661761 & 2 & 1 & (0;2^6,3^8) \\
905177 & 8 & 1 & (2;2^4,3^4) &	2782261 & 3 & 1 & (2;2^{14},3^4) \\
966125 & 5 & 1 & (1;3^4,5^2) &	2803712 & 2 & 1 & (2;3^6,4) \\
1075648 & 7 & 1 & (2;2^9,14) &	2847089 & 2 & 1 & (1;2^6,3^6) \\
1081856 & 7 & 1 & (2;2^8,4^2) &	2936696 & 2 & 1 & (1;2^4,3^8) \\
1202933 & 5 & 1 & (2;3^6) &	3195392 & 2 & 1 & (2;3^{10}) \\
1229312 & 7 & 1 & (2;2^8,4^2,7^2) &	3319769 & 2 & 1 & (1;2^{10},3^6) \\
1241125 & 5 & 1 & (2;3^4,5^2) &	3389609 & 2 & 1 & (1;2^{10},3^6) \\
1259712 & 3 & 1 & (0;2^9,18) &	3697873 & 2 & 1 & (1;2^{10},3^8) \\
1312625 & 4 & 1 & (1;2^4,5^4) &	4125937 & 2 & 1 & (2;2^{10},3^8) \\
1387029 & 3 & 1 & (1;2^8) &	4254689 & 2 & 1 & (2;2^8,3^{10})
\end{array}
\] \\
\end{footnotesize}
\textbf{Table 4.6}: Shimura curves with $[F:\Q]=6$
\end{center}

\begin{center}
\begin{footnotesize}
\[
\begin{array}{cccc|cccc}
d_F & D & N & \sigma & d_F & D & N & \sigma \\
\hline
20134393\T & 1 & 1 & (0;2^5,3^3) &	39829313 & 1 & 1 & (2;2^6,3^4) \\
25164057 & 1 & 1 & (0;2^5,3^5) &	41153941 & 1 & 1 & (0;2^{10},3^7) \\
25367689 & 1 & 1 & (0;2^7,3^3) &	41455873 & 1 & 1 & (1;2^{10},3^4) \\
28118369 & 1 & 1 & (0;2^7,3^4) &	41783473 & 1 & 1 & (1;2^{10},3^4) \\
30653489 & 1 & 1 & (1;2^5,3^4) &	42855577 & 1 & 1 & (1;2^9,3^5) \\
31056073 & 1 & 1 & (0;2^7,3^5) &	43242544 & 1 & 1 & (1;2^{10},3^5) \\
32354821 & 1 & 1 & (0;2^8,3^5) &	43723857 & 1 & 1 & (2;2^7,3^5) \\
32567681 & 1 & 1 & (0;2^9,3^4) &	46643776 & 1 & 1 & (2;2^{15},3^3) \\
34554953 & 1 & 1 & (1;2^6,3^4) &	52011969 & 1 & 1 & (2;2^7,3^9) \\
35269513 & 1 & 1 & (0;2^9,3^5) &	55073801 & 1 & 1 & (2;2^{11},3^6) \\
39610073 & 1 & 1 & (1;2^9,3^4) &	
\end{array}
\] \\
\end{footnotesize}
\textbf{Table 4.7}: Shimura curves with $[F:\Q]=7$
\end{center}

\renewcommand{\baselinestretch}{1}

\section{Comments}

In addition to employing the formula for the signature as in \S 2, the data was also verified as follows.  For all groups $\Gamma_0^{\frakD}(\frakN)$ with $N(\frakN) > 1$, we employed an algorithm to compute a fundamental domain for the group which also explicitly computes a presentation for the group \cite{VoightFuchsiangroups}, and in particular the signature.  This extra computation is independent and verified the correctness of the signature in each case.

We conclude by noting some points of interest in the data.  The first ambiguity in the tables occurs with $[F:\Q]=2$, $d_F=8$, $D=2$, and $N=49$.  One can choose $\frakN$ as either the square of an ideal of norm $7$, in which case the corresponding curve has signature $(2;3^4)$, or choose $\frakN$ as the ideal generated by $7$, in which case the signature is $(1;3^8)$.  Note that for $[F:\Q]=4$, $d_F=8468$, $D=2$, and $N=1$, we see that the signature even depends on which prime of norm $2$ one chooses for $\frakD$.

The largest elliptic cycle has order $18$, with $F \cong \Q(\zeta_{36})^+$ having $[F:\Q]=6$ and $d_F=108^3=1259712$.  The largest elliptic cycle of prime order is $13$.  The curve with $[F:\Q]=5$, $d_F=341692$, and $D=N=1$ has $20$ elliptic cycles of order $2$.  The largest level $D=232$ and discriminant $N=232$ occur with $F=\Q(\zeta_7)^+$ having $[F:\Q]=3$ and $d_F=49$.

Finally, we note that whenever $[F:\Q]$ is odd, $D=N=1$, and $g \geq 1$, the Jacobian of the curve $X_0^{(1)}(1)$ has everywhere good reduction.  For example, when $[F:\Q]=3$, this occurs for the fields with $d_F=961,1016,1101,\dots$.  

It is interesting to note the largest and smallest areas that occur for given genus.  It is well-known that the $(2,3,7)$-triangle group (of genus zero) has the smallest coarea for any arithmetic Fuchsian group.  For $g=1$, the curve over $\Q(\sqrt{13})$ with $D=4$ and $N=1$ has smallest area $1/2$; for $g=2$ the curve over $\Q$ with $D=26$ has smallest area $2$.  The groups with largest areas all occur with $[F:\Q]=4$, as follows.  For genus $0$, the curve with $d_F=7168$, $D=7$ and $N=1$ has area $A=17/2$.  For genus $1$, the curve with $d_F=30056$, $D=2$ and $N=1$ has area $38/3$.  Finally, for genus $2$, the curve with $d_F=2000$, $D=4$ and $N=25$ has area $15$.  These records may get larger if one allows an arbitrary arithmetic Fuchsian group: indeed, already the classical modular curve $X(5)$ has genus zero and area $10$.

Finally, it would be interesting to obtain equations for these curves, as well as to locate the elliptic (and parabolic) points.  For genus $0$, one could list the conic (or equivalently, the set of ramified places of the associated quaternion algebra); for genus $1$, one could provide a Weierstrass equation for the Jacobian of the curve, or in the case of low index an equation for the curve itself; and for a genus $2$ curve, one could give a hyperelliptic equation for the curve.


\begin{thebibliography}{999}

\bibitem{AB}
Montserrat Alsina and Pilar Bayer, \emph{Quaternion orders, quadratic forms, and Shimura curves}, CRM Monograph Series, vol.~22, Amer.\ Math.\ Soc., Providence, 2004. 

\bibitem{CLY}
Kok Seng Chua, Mong Lung Lang, and Yifan Yang, \emph{On Rademacher's conjecture: congruence subgroups of genus zero of the modular group}, J.~Algebra \textbf{277} (2004), no.~1, 408--428.

\bibitem{Cohen2}
Henri Cohen, \emph{Advanced topics in computational number theory}, Grad.\ Texts in Math., vol.~193, Springer-Verlag, New York, 2000.

\bibitem{CDO}
Henri Cohen, Francisco Diaz y Diaz, and Michel Olivier, \emph{Constructing complete tables of quartic fields using Kummer theory}, Math.\ Comp.\ \textbf{72} (2003), no.~242, 941--951.

\bibitem{CWZ}
J\'anos A.~Csirik, Joseph L.~Wetherell, and Michael E.~Zieve, On the genera of $X_0(N)$, to appear in J.\ Number Theory, 
\texttt{math.NT/0006096}.

\bibitem{DH}
Fred Diamond and Jerry Shurman, \emph{A first course in modular forms}, Grad. Texts in Math., vol.~228, Springer-Verlag, New York, 2005.

\bibitem{Dokchitser}
Tim Dokchitser, \emph{Computing special values of motivic $L$-functions}, Experiment.\ Math.\ \textbf{13} (2004), no.~2, 137--149. 

\bibitem{ElkiesSCC}
Noam D.~Elkies, \emph{Shimura curve computations}, Algorithmic number theory (Portland, OR, 1998), Lecture notes in Comput.\ Sci., vol. 1423, Springer, Berlin, 1998, 1--47.

\bibitem{JL}
Herv\'e Jacquet and Robert Langlands, \emph{Automorphic forms on $\GL(2)$}, Lecture Notes in Math., vol.~114, Springer-Verlag, Berlin-New York, 1970.

\bibitem{Johansson}
Stefan Johansson, \emph{Genera of arithmetic Fuchsian groups}, Acta Arith.\ \textbf{86} (1998), no.~2, 171--191. 

\bibitem{Katok}
Svetlana Katok, \emph{Fuchsian groups}, University of Chicago Press, Chicago, 1992.

\bibitem{LMR}
Darren D.~Long, Colin Maclachlan, Alan Reid, Arithmetic Fuchsian groups of genus zero, \emph{Pure Appl.~Math.~Q.}\ \textbf{2} (2006), no.~2, 569--599.

\bibitem{Louboutin}
St\'ephane Louboutin, Computation of relative class numbers of CM-fields, \emph{Math. Comp.}\ \textbf{66} (1997), no.~219, 1185--1194.

\bibitem{MR}
Colin Maclachlan and Gerhard Rosenberger, \emph{Two-generator arithmetic Fuchsian groups. II.}, Math.\ Proc.\ Cambridge Philos.\ Soc.\ \textbf{111} (1992), no.~1, 7--24.

\bibitem{Martinet}
Jacques Martinet, Petits discriminants des corps de nombres, \emph{Number theory days (Exeter, 1980)}, London Math. Soc. Lecture Note Ser., vol.~56, Cambridge Univ. Press, Cambridge-New York, 1982, 151--193.

\bibitem{Odlyzko}
Andrew M.~Odlyzko, Bounds for discriminants and related estimates for class numbers, regulators and zeros of zeta functions: a survey of recent results, \emph{S\'em. Th\'eor. Nombres Bordeaux} (2) \textbf{2} (1990), no.~1, 119--141. 

\bibitem{Sarnak}
Peter Sarnak, Notes on the Generalized Ramanujan Conjectures, \emph{Clay Math Proceedings}, vol.~4, 2005, 659--685.

\bibitem{Schneider}
Volker Schneider, Die elliptischen Fixpunkte zu Modulgruppen in Quaternionenschief\-k\"orpern, \emph{Math.~Ann.} \textbf{217} (1975), no.~1, 29--45.

\bibitem{Selberg}
Atle Selberg, On the estimation of Fourier coefficients of modular forms, \emph{Proc.\ Sympos.\ Pure Math.}, Vol.~VIII, Amer.\ Math.\ Soc., Providence, 1965, 1--15.

\bibitem{Shimizu}
Hideo Shimizu, On zeta functions of quaternion algebras, \emph{Ann.~of Math.~(2)} \textbf{81}, 1965, 166--193.

\bibitem{MShimura}
Mahoro Shimura, Defining equations of modular curves $X\sb 0(N)$, \emph{Tokyo J.\ Math.} \textbf{18} (1995), no.\ 2, 443--456.

\bibitem{Takeuchi}
Kisao Takeuchi, Arithmetic Fuchsian groups with signature $(1;e)$, 
\emph{J.~Math.~Soc.~Japan}, \textbf{35} (1983), no.~3, 381--407. 

\bibitem{Thompson}
John G.~Thompson, A finiteness theorem for subgroups of $\PSL(2,\R)$ which are commensurable with $\PSL(2,\Z)$, \emph{The Santa Cruz conference on finite groups (Univ. California, Santa Cruz, Calif., 1979)}, Proc.\ Sympos.\ Pure Math., vol.~37, Amer.\ Math.\ Soc., Providence, 1980, 533--555.

\bibitem{Vigneras}
Marie-France Vign\'eras, \emph{Arithm\'etique des alg\`ebres de quaternions}, Lecture Notes in Math., vol.~800, Springer, Berlin, 1980. 

\bibitem{Vigneras1-4}
Marie-France Vign\'eras, \emph{Quelques remarques sur la conjecture $\lambda_1 \geq 1/4$}, Seminar on number theory, Paris 1981--82, Progr.~Math.\ \textbf{38} (1983) 321--343.

\bibitem{Voight}
John Voight, \emph{Enumeration of totally real number fields of bounded root discriminant}, to appear in ANTS VIII.

\bibitem{Voightonline}
John Voight, Shimura curves of genus at most two, \\ \texttt{http://www.cems.uvm.edu/\~{}voight/shim-tables/}\ .

\bibitem{VoightFuchsiangroups}
John Voight, \emph{Computing fundamental domains for Fuchsian groups}, submitted to J.\ Th\'eorie de Nombres
de Bordeaux.

\bibitem{Zograf}
Peter Zograf, A spectral proof of Rademacher's conjecture for congruence subgroups of the modular group, \emph{J.\ reine angew.\ Math.} \textbf{414} (1991), 113--116.

\end{thebibliography}
\end{document}